\documentstyle[12pt]{article}

\newcommand{\h}[2]{h^{#1}_{#2}}		
\newcommand{\Qh}[1]{Q_h^{#1}} 		
\newcommand{\al}[2]{\alpha^{#1}_{#2}}	
\newcommand{\Qal}[1]{Q_\alpha^{#1}} 	
\newcommand{\QLa}[1]{Q_\Lambda^{#1}} 	

\renewcommand\a{\alpha}
\newcommand\la{\lambda}
\def\ddelta{\,\delta}
\newcommand\La{\Lambda}
\newcommand\bZ{{\bf Z}}
\newcommand\bC{{\bf C}}
\newcommand\cF{{\cal F}}
\newcommand\cW{{\cal W}}

\newcommand\rb{\sqrt{\beta}}
\newcommand\rbi{-{1\over\sqrt{\beta}}}

\newcommand\hb{\hfill\break}

\renewcommand\/{\over}
\renewcommand\({\left(}
\renewcommand\]{\,\right]}
\renewcommand\[{\left[\,}
\renewcommand\){\right)}
\def\<{\langle}
\renewcommand\>{\rangle}
\def\qed{\hfill\fbox{}}
\def\proof{\noindent{\it Proof.\quad}}

\newcommand{\Exp}[1]{\exp\left\{#1\right\}}
\renewcommand{\:}{{\textstyle{\circ\atop\circ}}}
\newcommand{\qVir}{${\cal V}ir\hspace{-.03in}_{q,t}\,$}
\newcommand{\qW}{$q$-${\cal W}_N$}

\newcommand{\be}{\begin{equation}}
\newcommand{\ee}{\end{equation}}
\newcommand{\bea}{\begin{eqnarray}}
\newcommand{\eea}{\end{eqnarray}}
\renewcommand{\&}{&\!\!\!\!\!\!\!\! &}

\newcommand{\eq}[1]{(\ref{#1})}

\def\boxline#1{\vbox{\hrule\hbox{\vrule\vbox{#1}\vrule}\hrule}}
\def\boxNW#1{\vbox{\hrule\hbox{\vrule\vbox{#1}}}}
\def\boxES#1{\vbox{\hbox{\vbox{#1}\vrule}\hrule}}
\def\Between#1#2#3#4{ 
\raise-#1mm\vbox to#1mm{\hsize #2mm \vbox{\vskip #3mm\centerline{#4} } }
}
\def\Square#1#2#3#4{ 
\raise-#1mm\boxline{
\vbox to#1mm{\hsize #1mm \vbox{\vskip #2mm\noindent\hskip4pt {#3}
			       \vskip-#2mm\vskip-8pt\centerline{#4} }} }
}
\def\Young#1#2#3#4#5#6#7#8#9{
\raise-#9mm\boxNW{\vbox to#1mm{\hsize#6mm
	\vbox{\vskip#8mm\noindent\hskip#8mm$\;\lambda$} }}
\kern-#6mm
\raise-#1mm\boxES{\vbox to #5mm{\hsize #7mm $ $}}\kern-.4pt
\raise-#2mm\boxES{\vbox to #5mm{\hsize #7mm $ $}}\kern-.4pt
\raise-#3mm\boxES{\vbox to #5mm{\hsize #7mm $ $}}\kern-.4pt
\raise-#4mm\boxES{\vbox to #5mm{\hsize #7mm $ $}}\kern-.4pt
\raise-#5mm\boxES{\vbox to #5mm{\hsize #7mm $ $}}
}
\def\Galilei{
  \Young{10}8642{15}32{9.9}
\Between{10}{20}3{$\longmapsto$}
 \Square{14}7{$r$}{$s$} \kern-.4pt  
  \Young{10}8642{15}32{9.9}
\Between{10}{10}{10}{\quad .}
}
\def\generalYoung{
\vskip.25cm
\noindent
\makebox[  3cm]{ }
\makebox[  2cm]{$s_1$}\hskip-.4pt
\makebox[1.7cm]{$s_2$}
\makebox[1.4cm]{ }
\makebox[1.4cm]{$s_{N-2}$}\hskip-.35pt
\makebox[1.3cm]{$s_{N-1}$}
\hfill\break
 \makebox[  3cm][r]{$\hfill\lambda=$}
\framebox[  2cm][l]{\rule[  -1cm]{0cm}{  2cm}$r_1$}\hskip-.4pt
\framebox[1.7cm][l]{\rule[-0.7cm]{0cm}{1.7cm}$r_2$}
 \makebox[1.4cm]			    {\raisebox{.25cm}{$\cdots\cdots$}}
\framebox[1.4cm][l]{\rule[-0.4cm]{0cm}{1.4cm}\raisebox{.25cm}{$r_{N-2}$}
							     }\hskip-.4pt
\framebox[1.3cm][l]{\rule[-0.2cm]{0cm}{1.2cm}\raisebox{.25cm}{$r_{N-1}$}}
\makebox[1cm][r]{.}
\vskip.3cm
}

\begin{document}


\renewcommand{\thefootnote}{\fnsymbol{footnote}}
\font\csc=cmcsc10 scaled\magstep1

{\baselineskip=14pt
 \rightline{
 \vbox{
       \hbox{YITP/U-95-34}
       \hbox{DPSU-95-9}
       \hbox{UT-718}
       \hbox{August 1995}
}}}

\vskip 11mm
\begin{center}

{\large\bf Quantum $\cW_N$ Algebras and Macdonald Polynomials}

\vspace{15mm}

{\csc Hidetoshi AWATA},\footnote{JSPS fellow}$^{1}$
{\csc Harunobu KUBO},$^{*2}$
{\csc Satoru ODAKE}$^3$
\\ and
{\csc Jun'ichi SHIRAISHI}$^4$
{\baselineskip=15pt
\it\vskip.35in 
\setcounter{footnote}{0}\renewcommand{\thefootnote}{\arabic{footnote}}
\footnote{e-mail address : awata@yisun1.yukawa.kyoto-u.ac.jp}
Uji Research Center, Yukawa Institute for Theoretical Physics \\
Kyoto University, Uji 611, Japan
\vskip.1in 
\footnote{e-mail address : kubo@danjuro.phys.s.u-tokyo.ac.jp}
Department of Physics, Faculty of Science \\
University of Tokyo, Tokyo 113, Japan \\
\vskip.1in 
\footnote{e-mail address : odake@yukawa.kyoto-u.ac.jp}
Department of Physics, Faculty of Science \\
Shinshu University, Matsumoto 390, Japan\\
\vskip.1in 
\footnote{e-mail address : shiraish@momo.issp.u-tokyo.ac.jp}
Institute for Solid State Physics, \\
University of Tokyo, Tokyo 106, Japan \\
}
\end{center}

\vspace{7mm}
\begin{abstract}

We derive a quantum deformation of the $\cW_N$ algebra and
its quantum Miura transformation,
whose singular vectors realize the Macdonald polynomials.

\end{abstract}

\vspace{10mm}
q-alg/9508011

\vfill\eject
\setcounter{footnote}{0}\renewcommand{\thefootnote}{\arabic{footnote}}


\section{Introduction}


The excited states of the Calogero-Sutherland model \cite{rSu}
and its relativistic model
(the trigonometric limit of the Ruijsenaars model) \cite{rR}
are described by the Jack polynomials \cite{rSt}
and their $q$-analog (the Macdonald polynomials) \cite{rM}, respectively.
Since the Jack polynomials coincide with
certain correlation functions of $\cW_N$ algebra \cite{rMY,rAMOS},
it is natural to expect that the Macdonald polynomials
are also realized by those of a deformation of $\cW_N$ algebra.

In a previous paper \cite{rSKAO},
we derived a quantum Virasoro algebra whose singular vectors are 
some special kinds of Macdonald polynomials.
On the other hand,
E.~Frenkel and N.~Reshetikhin succeeded in constructing
the Poisson $\cW_N$ algebra
and its quantum Miura transformation
in the analysis of the $U_q(\widehat{sl_N})$ algebra at the critical level
\cite{rFR}.
Like the classical case \cite{rFL},
these two works, $q$-Virasoro and $q$-Miura transformation,
are essential to find and study a quantum $\cW_N$ algebra.
In this article,
we present a {\qW} algebra
whose singular vectors realize the general Macdonald polynomials.

This paper is arranged as follows:
In section 2, we define a quantum deformation of $\cW_N$ algebras
and its quantum Miura transformation.
The screening currents and a vertex operator are derived
in section 3 and 4.
A relation with the Macdonald polynomials is obtained in section 5.
Section 6 is devoted to conclusion and discussion.
Finally we recapitulate the $q$-Virasoro algebra and
the integral formula for the Macdonald polynomials in appendices.


\section{Quantum deformation of $\cW_N$ algebra}


We start with defining
a new quantum deformation of the $\cW_N$ algebra
by quantum Miura transformation.


\subsection{Quantum Miura transformation}


\newcommand\frenkel{
We found this commutation relation by comparing
the Poisson bracket in Frenkel-Reshetikhin's work \cite{rFR}
and the commutator in ours \cite{rSKAO}.
The oscillator $a_n$ used in \cite{rSKAO} is given by
$a_n = -n h^1_n p^{-n/2}/(1-t^n)$ and
$a_{-n} = n h^1_{-n} p^{n/2}(1+p^n)/(1-t^{-n})$ for $n>0$.
}

First we define fundamental bosons ${\h in}$ and ${\Qh i}$
for $i=1,2,\cdots,N$ and  $n\in\bZ$ such that\footnote{\frenkel}
\bea
[{\h in},{\h jm}] &\!\!=\!\!& -{1\/n}(1-q^n)(1-t^{-n})
	{1-p^{(\ddelta_{ij}N-1)n}\/1-p^{Nn}} p^{Nn\theta(i<j)}
	\ddelta_{n+m,0},\cr
[{\h i0},{\Qh j}] &\!\!=\!\!& \ddelta_{ij} - {1\/N},\qquad\quad
\sum_{i=1}^N p^{in} {\h in} = 0,\qquad\quad
\sum_{i=1}^N {\Qh i} =0,
\eea
with $q$, $t\equiv q^\beta\in\bC$ and $p\equiv q/t$.
Here $\theta(P)\equiv 1$ or $0$
if the proposition $P$ is true or false, respectively.
This bosons correspond to the weights of the vector representation
$h_i$ whose inner-product is $(h_i\cdot h_j)=(\ddelta_{ij}N-1)/N$.


Let us define fundamental vertices $\La_i(z)$ and
{\qW} generators $W^i(z)$ for $i=1,2,\cdots,N$ as follows:
\bea
\La_i(z) &\!\!\equiv\!\!& :\Exp{ \sum_{n\neq 0}{\h in} z^{-n} }:
	q^{\rb{\h i0}} p^{{N+1\/2}-i},\cr
W^i(zp^{1-i\/2}) &\!\!\equiv\!\!& \sum_{1\leq j_1<\cdots<j_i\leq N}
	:\La_{j_1}(z) \La_{j_2}(zp^{-1}) \cdots \La_{j_i}(zp^{1-i}):,
\eea
and $W^0(z)\equiv 1$.
Here $:*:$ stands for the usual bosonic normal ordering such that
the bosons ${\h in}$ with non-negative mode $n\geq 0$ are in the right. 
Note that
\be
W^N(zp^{1-N\/2}) =
	\,:\!\La_1(z) \La_2(zp^{-1}) \cdots \La_N(zp^{1-N})\!:\, = 1.
\ee
If we take the limit $t\rightarrow 1$ with $q$ fixed,
the above generators reduce to those of Ref.\ \cite{rFR}.
These generators are obtained by the following quantum Miura transformation:
\be
:\!\(p^{D_z} - \La_1(z)\) \(p^{D_z} - \La_2(zp^{-1})\) \cdots
\(p^{D_z} - \La_N(zp^{1-N})\)\!:\,
= \sum_{i=0}^N (-1)^i W^i(zp^{1-i\/2}) p^{(N-i)D_z},
\label{e:qMiura}
\ee
with $D_z \equiv z{\partial\/\partial z}$.
Remark that $p^{D_z}$ is the $p$-shift operator such that
$p^{D_z} f(z) = f(pz)$.


\subsection{Relations of {\qW} generators}


Next we give the algebra of the above {\qW} generators.
Let $W^i(z) = \sum_{n\in\bZ} W^i_n z^{-n}$.
Let us define a new normal ordering $\:*\:$
for the {\qW} generators as follows:
\bea
&\!\!&\!\!\!\!\!\!
  \: W^i(rw)W^j(w)\: \cr
&\!\!\equiv\!\!&
  \oint{dz\/2\pi iz}\left\{
 {     1\/1-{rw/z}}f^{ij}\({w\/z}\)W^i(z)W^j(w)
+{{z/rw}\/1-{z/rw}}                W^j(w)W^i(z)f^{ji}\({z\/w}\)
\right\}\cr
&\!\!=\!\!&
\sum_{n\in\bZ}\sum_{m\geq 0}\sum_{\ell=0}^m f^{ij}_\ell \left\{
  r^{  m-\ell}\cdot W^i_{-m}   W^j_{n+m}
 +r^{\ell-m-1}\cdot W^j_{n-m-1}W^i_{m+1}
  \right\}w^{-n},
\eea
with
\bea
f^{ij}(x) &\!\!\equiv\!\!& \Exp{ \sum_{n>0}{1\/n}(1-q^n)(1-t^{-n})
	{1-p^{in}\/1-p^n}{1-p^{(N-j)n}\/1-p^{Nn}} p^{{j-i\/2}n} x^n },\cr
f^{ji}(x) &\!\!\equiv\!\!& f^{ij}(x),\qquad (i\leq j),
\eea
and
$f^{ij}(x)\equiv\sum_{\ell\geq 0}f^{ij}_\ell x^\ell$.
Here $(1-x)^{-1}$ stands for $\sum_{n\geq 0}x^n$.
Remark that this normal ordering $\:*\:$ is a generalization
of the following usual one $(*)$ used in the conformal field theory:
\bea
\(AB\)(w) &\!\!\equiv\!\!& \oint_w {dz\/2\pi i} {1\/z-w} A(z) B(w)\cr
&\!\!\equiv\!\!& \oint_0 {dz\/2\pi iz}
\left\{{1\/1-{w/z}}A(z)B(w) + {{z/w}\/1-{z/w}}B(w)A(z)\right\}.
\eea


The relation of the {\qW} generators should be
written in this normal ordering.
Here we present some examples of them.
The relation of $W^1(z)$ and $W^j(z)$ for $j\geq 1$ is
\bea
&\!\!\!\!\!\!&\!\!\!\!\!\!
f^{1j}\({w\/z}\)W^1(z)W^j(w) - W^j(w)W^1(z)f^{j1}\({z\/w}\) \\
&\!\!=\!\!&
  -{(1-q)(1-t^{-1})\/1-p}\left\{
  \ddelta\(p^{j+1\/2}{w\/z}\)W^{j+1}\(p^{1\/2}w\)
  -\ddelta\(p^{-{j+1\/2}}{w\/z}\)W^{j+1}\(p^{-{1\/2}}w\)
  \right\},\nonumber
\eea
with $\ddelta(x)\equiv \sum_{n\in\bZ} x^n$;
and that of $W^2(z)$ and $W^j(z)$ for $j\geq 2$ is
\bea
&\!\!\!\!\!\!&\!\!\!\!\!\!
f^{2j}\({w\/z}\)W^2(z)W^j(w) - W^j(w)W^2(z)f^{j2}\({z\/w}\) \\
&\!\!=\!\!&
  -{(1-q)(1-t^{-1})\/1-p}{(1-qp)(1-t^{-1}p)\/(1-p)(1-p^2)}
  \left\{\ddelta\(p^{ {j\/2}+1}{w\/z}\)W^{j+2}(p     w) \right. \cr
&\!\!\!\!\!\!&\hspace{66mm}
  \left.-\ddelta\(p^{-{j\/2}-1}{w\/z}\)W^{j+2}(p^{-1}w) \right\}\cr
&\!\!\!\!\!\!&
  -{(1-q)(1-t^{-1})\/1-p}
  \left\{\ddelta\(p^{{j\/2}}{w\/z}\)
         \: W^1(p^{-{1\/2}}z) W^{j+1}(p^{ {1\/2}} w)\: \right. \cr
&\!\!\!\!\!\!&\hspace{32mm}
  \left.-\ddelta\(p^{-{j\/2}}{w\/z}\)
         \: W^1(p^{ {1\/2}}z) W^{j+1}(p^{-{1\/2}} w)\: \right\}\cr
&\!\!\!\!\!\!&
  +{(1-q)^2(1-t^{-1})^2\/(1-p)^2}
  \left\{\ddelta\(p^{{j\/2}}{w\/z}\)
         \({p^2\/1-p^2}W^{j+2}(pw)+{1\/1-p^j} W^{j+2}(      w)\)\right.\cr
&\!\!\!\!\!\!&\hspace{35mm}
  \left.-\ddelta\(p^{-{j\/2}}{w\/z}\)
         \({p^j\/1-p^j}W^{j+2}( w)+{1\/1-p^2} W^{j+2}(p^{-1}w)\)\right\},
\nonumber
\eea
with $W^i(z) \equiv 0$ for $i>N$.
The main terms of
$$
f^{ij}\({w\/z}\)W^i(z)W^j(w) - W^j(w)W^i(z)f^{ji}\({z\/w}\)\qquad
(i\leq j)
$$
is
\bea
&\!\!\!\!\!\!&
  -{(1-q)(1-t^{-1})\/1-p}\sum_{k=1}^{{\rm min} \(i,N-j\)}
\prod_{\ell=1}^{k-1}{(1-qp^\ell)(1-t^{-1}p^\ell)\/(1-p^\ell)(1-p^{\ell+1})}\cr
&\!\!\!\!\!\!&\!\!\!\!\!\! \times
\left\{
 \ddelta\(p^{{j-i\/2}+k}{w\/z}\)\:W^{i-k}(p^{-{k\/2}}z)W^{j+k}(p^{ {k\/2}}w)\:
-\ddelta\(p^{{i-j\/2}-k}{w\/z}\)\:W^{i-k}(p^{ {k\/2}}z)W^{j+k}(p^{-{k\/2}}w)\:
\right\}.
\nonumber
\eea


\newcommand\AC{
In these kinds of formulae we use
$ \Exp{-\sum_{n>0}x^n   /n} = 1-x = 
-x\Exp{-\sum_{n>0}x^{-n}/n}$.}

To obtain the above relations, the fundamental formula is
\bea
f^{11}\({w\/z}\)\La_i(z)\La_i(w) &\!\!-\!\!& \La_i(w)\La_i(z)f^{11}\({z\/w}\)
= 0,\cr
f^{11}\({w\/z}\)\La_i(z)\La_j(w) &\!\!-\!\!& \La_j(w)\La_i(z)f^{11}\({z\/w}\)
\cr
&\!\!=\!\!&{(1-q)(1-t^{-1})\/1-p}\(\ddelta\({w\/z}\)-\ddelta\(p{w\/z}\)\)
:\La_i(z)\La_j(w):,
\nonumber
\eea
for $i<j$, here we use\footnote{\AC}
\bea
&\!\!\!\!\!\!&\!\!\!\!\!\!
\Exp{ \sum_{n>0}{1\/n}(1-q^{ n})(1-t^{-n})x^{ n} }-
\Exp{ \sum_{n>0}{1\/n}(1-q^{-n})(1-t^{ n})x^{-n} }\cr
&\!\!\!\!\!\!&\!\!\!\!\!\! \hspace{55mm}
= {(1-q)(1-t^{-1})\/1-p}\(\ddelta(x)-\ddelta(px)\).
\eea
To calculate the general relations, the following formulae are useful:
\bea
&\!\!\!\!\!\!&\!\!\!\!\!\!
\Exp{ \sum_{n>0}{1\/n}(1-q^{ n})(1-t^{-n})(1+r^{ n})x^{ n} }-
\Exp{ \sum_{n>0}{1\/n}(1-q^{-n})(1-t^{ n})(1+r^{-n})x^{-n} }\cr
&\!\!\!\!\!\!&\!\!\!\!\!\! \hspace{30mm}
={(1-q)(1-t^{-1})\/(1-p)(1-r)}
\left\{
 (1-qr)(1-t^{-1}r) {\ddelta( x)-\ddelta(prx)\/1-pr} \right.\cr
&\!\!\!\!\!\!&\!\!\!\!\!\! \hspace{70mm}\left.
-(r-q )(r-t^{-1} ) {\ddelta(rx)-\ddelta(p x)\/r-p }
\right\},
\label{e:formula.3}
\eea
with $r\neq 0$;
For $r=1$ or $p^{\pm1}$,
the right hand side of \eq{e:formula.3} should be understood as the limit
$r\rightarrow 1$ or $p^{\pm1}$, respectively;
And
$f^{ij}(x) = \prod_{k=1}^i f^{1j}(p^{{i+1\/2}-k}x)$
for $i\leq j$.


\subsection{Example of $q$-$\cW_3$}


$N=2$ case is {\qVir} studied in Ref. \cite{rSKAO} (see appendix A).\\
Here we give an example when $N=3$.
The generators are
\bea
W^1(z) &\!\!=\!\!& \La_1(z) + \La_2(z) + \La_3(z),\cr
W^2(z) &\!\!=\!\!&
\La_1(zp^{1\/2})\La_2(zp^{-{1\/2}}) +
\La_1(zp^{1\/2})\La_3(zp^{-{1\/2}}) +
\La_2(zp^{1\/2})\La_3(zp^{-{1\/2}}).
\eea
The relation of these generators is
\bea
\&\hspace{-5mm}
f^{11}\({w\/z}\) W^1(z) W^1(w) - W^1(w) W^1(z) f^{11}\({z\/w}\) \cr
\&\hspace{10mm}
= -{(1-q)(1-t^{-1})\/1-p}\left\{
\ddelta\({w\/z}p     \) W^2\(wp^{  1\/2 }\)-
\ddelta\({w\/z}p^{-1}\) W^2\(wp^{-{1\/2}}\)\right\},\cr
\&\hspace{-5mm}
f^{12}\({w\/z}\) W^1(z) W^2(w) - W^2(w) W^1(z) f^{21}\({z\/w}\) \cr
\&\hspace{10mm}
= -{(1-q)(1-t^{-1})\/1-p}\left\{
\ddelta\({w\/z}p^{ {3\/2}}\)-
\ddelta\({w\/z}p^{-{3\/2}}\)\right\},\cr
\&\hspace{-5mm}
f^{22}\({w\/z}\) W^2(z) W^2(w) - W^2(w) W^2(z) f^{22}\({z\/w}\) \cr
\&\hspace{10mm}
= -{(1-q)(1-t^{-1})\/1-p}\left\{
\ddelta\({w\/z}p     \) W^1\(zp^{-{1\/2}}\)-
\ddelta\({w\/z}p^{-1}\) W^1\(zp^{  1\/2 }\)\right\},\nonumber
\eea
with
\bea
f^{11}(x) &\!\!=\!\!&
\Exp{ \sum{1\/n}(1-q^n)(1-t^{-n}){1-p^{2n}\/1-p^{3n}}x^n }
=f^{22}(x),\cr
f^{12}(x) &\!\!=\!\!&
\Exp{ \sum{1\/n}(1-q^n)(1-t^{-n}){1-p^n\/1-p^{3n}}p^{n\/2}x^n }
=f^{21}(x).\nonumber
\eea
Note that there is no difference between $W^1$ and $W^2$
in algebraically.


\subsection{Highest weight module of {\qW} algebra}


Here we refer to the representation of the {\qW} algebra.
Let $|\la\>$ be the highest weight vector of the {\qW} algebra
which satisfies
$W^i_n|\la\>=0$ for $n>0$ and $i=1,2,\cdots,N-1$ and
$W^i_0|\la\>= \la^i |\la\>$ with $\la^i\in\bC$.
Let $M_\la$ be the Verma module over the {\qW} algebra
generated by $|\la\>$.
The dual module $M_\la^*$ is generated by $\<\la|$ such that
$\<\la|W^i_n=0$ for $n<0$ and $\<\la|W^i_0= \la^i\<\la|$.
The bilinear form $M_\la^*\otimes M_\la\rightarrow\bC$
is uniquely defined by $\<\la|\la\>=1$.

A singular vector $|\chi\>\in M_\la$ is defined by
$W^i_n|\chi\>=0$ for $n>0$ and
$W^i_0|\chi\>= (\la^i+N^i) |\chi\>$
with $N^i\in\bC$.


\section{Screening currents and singular vectors}


Next we turn to the screening currents, a commutant of the {\qW} algebra,
which construct the singular vectors.


\subsection{Screening currents}



Let us introduce root bosons
${\al in} \equiv {\h in}-{\h {i+1}n}$ and
${\Qal i} \equiv {\Qh i}-{\Qh {i+1}}$ for $i=1,2,\cdots,N-1$.
Then they satisfies
\bea
[{\al in},{\al jm}] &\!\!=\!\!& -{1\/n}(1-q^n)(1-t^{-n})
\left\{(1+p^{-n})\ddelta_{i,j} - \ddelta_{i+1,j} -
p^{-n}\ddelta_{i-1,j}\right\}
	\ddelta_{n+m,0},\cr
[{\al i0},{\Qal j}] &\!\!=\!\!& 2\ddelta_{i,j}-\ddelta_{i+1,j}-\ddelta_{i-1,j},
\eea
and
\bea
[{\h in},{\al jm}] &\!\!=\!\!& {1\/n}(1-q^{-n})(1-t^{-n})
	\left\{ q^n\ddelta_{i,j} - t^n\ddelta_{i,j+1}\right\}
	\ddelta_{n+m,0},\cr
[{\h i0},{\Qal j}] &\!\!=\!\!& \ddelta_{i,j}-\ddelta_{i,j+1},\qquad
[{\al i0},{\Qh j}]  =  \ddelta_{i,j}-\ddelta_{i+1,j}.
\eea
Note that
$[{\h in} + p^n {\h {i+1}n}, {\al im}] = 0$.


By using these root bosons, we define screening currents as follows:
\bea
S^i_+(z) &\!\!\equiv\!\!&
:\Exp{ \sum_{n\neq0}{{\al in}\/1-q^n} z^{-n} }:
e^{\rb{\Qal i}} z^{\rb{\al i0}},\cr
S^i_-(z) &\!\!\equiv\!\!&
:\Exp{ -\sum_{n\neq0}{{\al in}\/1-t^n} z^{-n} }:
e^{\rbi{\Qal i}} z^{\rbi{\al i0}}.
\eea
Then we have

\proclaim Proposition.
The screening currents satisfy
\bea
\&\hspace{-7mm}
\[\,:\!\(p^{D_z} - \La_1(z)\)\(p^{D_z} - \La_2(zp^{-1})\)\cdots
      \(p^{D_z} - \La_N(zp^{1-N})\)\!:\,,
S^i_\pm(w)\]\cr
\&
= (1-q^{\pm1})(1-t^{\mp1}){d\/d_{q\atop t}w}
\,:\!\(p^{D_z} - \La_1(z)\)\cdots\(p^{D_z} - \La_{i-1}(zp^{2-i})\) \cr
\&\hspace{5mm}\times
w \ddelta\({w\/z}p^{i-1}\) A^i_\pm(w) p^{D_z}
\(p^{D_z} - \La_{i+2}(zp^{-1-i})\)\cdots\(p^{D_z} - \La_N(zp^{1-N})\) \!:\,,
\nonumber
\eea
with
\bea
A^i_+(w) &\!\!=\!\!&
:\Exp{ \sum_{n\neq0}{{\h in}-q^n{\h {i+1}n}\/1-q^n}w^{-n} }:\,
e^{\rb{\Qal i}} w^{\rb{\al i0}} q^{\rb{\h{i+1}0}} p^{{N+1\/2}-i-1},\cr
A^i_-(w) &\!\!=\!\!&
:\Exp{ -\sum_{n\neq0}{t^n{\h in}-{\h {i+1}n}\/1-t^n}w^{-n} }:\,
e^{\rbi{\Qal i}} w^{\rbi{\al i0}} q^{\rb{\h i0}} p^{{N+1\/2}-i}.\nonumber
\eea

\noindent
Here ${d\/d_\xi w} f(w) \equiv (f(w)-f(\xi w))/((1-\xi)w)$.

\proof
First, we have
\bea
[\La_i(z),S^j_+(w)] &\!\!=\!\!&
(t-1)\ddelta_{i,j}\ddelta\({w\/z}q\):\La_j(z) S^j_+(w):\cr
&&+
(t^{-1}-1)\ddelta_{i,j+1}\ddelta\({w\/z}\):\La_{j+1}(z) S^j_+(w):,\cr
[\La_i(z),S^j_-(w)] &\!\!=\!\!&
(q^{-1}-1)\ddelta_{i,j}\ddelta\({w\/z}\):\La_j(z) S^j_-(w):\cr
&&+
(q-1)\ddelta_{i,j+1}\ddelta\({w\/z}t\):\La_{j+1}(z) S^j_-(w):.
\eea
Here we use the following formula:
\be
q^{\mp1}
\Exp{ \pm\sum_{n>0}{1\/n}(1-q^n   )x^n    } -
\Exp{ \pm\sum_{n>0}{1\/n}(1-q^{-n})x^{-n} } =
(q^{\mp1}-1)\ddelta\(xq^{1\mp1\/2}\).
\ee
The operator parts are
\bea
:\La_j(wq)     S^j_+(w): &\!\!=\!\!& A^j_+(wq) p,\qquad
:\La_{j+1}(w)  S^j_+(w):\,=  A^j_+(w),\cr
:\La_j(w)      S^j_-(w): &\!\!=\!\!& A^j_-(w),\qquad
:\La_{j+1}(wt) S^j_-(w):\,=  A^j_-(wt) p^{-1}.
\eea
Next,
\bea
[\La_i(z)+\La_{i+1}(z), S^i_\pm(w)] &\!\!=\!\!& -(1-q^{\pm1})(1-t^{\mp1})
{d\/d_{q\atop t}w}\left\{w\ddelta\({w\/z}\)A^i_\pm(w)\right\},\cr
[:\La_i(z)\La_{i+1}(zp^{-1}):, S^i_\pm(w)] &\!\!=\!\!& 0.
\eea
Hence, 
\bea
\&
\[\,:\!\(p^{D_z} - \La_i(z)\)\(p^{D_z} - \La_{i+1}(zp^{-1})\)\!:\,,
S^i_\pm(w)\]\cr
\&\hspace{35mm}
= (1-q^{\pm1})(1-t^{\mp1})
{d\/d_{q\atop t}w}\left\{w\ddelta\({w\/z}\)A^i_\pm(w)\right\} p^{D_z}.
\eea
This gives us the proposition. \qed

Therefore,
the screening currents $S^i_\pm(z)$
commute with any {\qW} generators up to total difference.
Thus we obtain

\proclaim Theorem.
Screening charges $\oint dz S^i_\pm(z)$
commute with any {\qW} generators.


\subsection{Singular vectors}


Let $\cF_\a$ be the boson Fock space
generated by the highest weight state $|\a\>$ such that
${\al in} |0\> = 0$ for $n\geq0$ and
$|\a\> \equiv \exp\{\sum_{i=1}^{N-1}\a^i{\QLa i}\} |0\>$
with ${\QLa i}\equiv\sum_{j=1}^i {\Qh j}$.
Note that ${\al i0}|\a\> = \a^i|\a\>$.
And this state $|\a\>$ is also the highest weight state of the {\qW} algebra.

We denote the negative mode part of $S^i_+(z)$ as
$(S^i_+(z))_- \equiv \Exp{ \sum_{n<0}{{\al in}\/1-q^n} z^{-n} }$.
Then we have

\proclaim Proposition.
For a set of non-negative integers $s_a$ and $r_a\geq r_{a+1}\geq0$,
($a=1,\cdots,N-1$),
let
\bea
          \a_{r,s}^{a} &\!\!=\!\!&\rb(1+r_a-r_{a-1})\rbi(1+s_a),\qquad
r_0 = 0,\cr
\widetilde\a_{r,s}^{a} &\!\!=\!\!&\rb(1-r_a+r_{a+1})\rbi(1+s_a),\qquad
r_N = 0.
\label{e:weightalpha}
\eea
Then the singular vectors $|\chi_{rs}^+\>\in\cF_{\a_{rs}^+}$
are realized by the screening currents as follows:
\bea
\&\hspace{-5mm}
|\chi_{r,s}\> =
\oint\prod_{a=1}^{N-1}\prod_{j=1}^{r_a} {dx^a_j}\cdot
S^1_+(x^1_1)\cdots S^1_+(x^1_{r_1}) \cdots
S^{N-1}_+(x^{N-1}_1) \cdots S^{N-1}_+(x^{N-1}_{r_{N-1}})
|\widetilde\a_{r,s}\>\cr
&\!\!=\!\!&
\oint\prod_{a=1}^{N-1} \prod_{j=1}^{r_a} {dx^a_j\/x^a_j}\cdot
\prod_{a=1}^{N-1} \Pi\(\overline{x^a},px^{a+1}\) \Delta(x^a) C(x^a)
\prod_{j=1}^{r_a} (x^a_j)^{-s_a} (S^a_+(x^a_j))_-\cdot|\a_{r,s}\>\cr
&&\label{e:singular}
\eea
with $x^N=0$, $\overline x = 1/x$ and
\bea
\Pi(x,y) &\!\!=\!\!& \prod_{ij}
\Exp{ \sum_{n>0}{1\/n}{1-t^n\/1-q^n} x_i^n y_j^n },\qquad
\Delta(x) =
\prod_{i\neq j}^r\Exp{ -\sum_{n>0}{1\/n}{1-t^n\/1-q^n}{x_j^n\/x_i^n}},\cr
C(x) &\!\!=\!\!&  \prod_{i<j}^r
\Exp{\sum_{n>0}{1\/n}{1-t^n\/1-q^n}\({x_i^n\/x_j^n}-p^n{x_j^n\/x_i^n}\)}
\prod_{i=1}^r x_i^{(r+1-2i)\beta}.
\label{e:Delta}
\eea

\proof
The operator product expansion of the screening currents is
\bea
S^a_+(x) S^a_+(y) &\!\!=\!\!&
\Exp{-\sum_{n>0}{1\/n}{1-t^n\/1-q^n}(1+p^n){y^n\/x^n}}x^{2\beta}
:S^a_+(x) S^a_+(y):,\cr
S^a_+(x) S^{a\pm1}_+(y) &\!\!=\!\!&
\Exp{\sum_{n>0}{1\/n}{1-t^n\/1-q^n}p^{{1\pm1\/2}n}{y^n\/x^n}}x^{-\beta}
:S^a_+(x) S^{a\pm1}_+(y):.
\eea
Since 
\bea
S^a_+(x_1)\cdots S^a_+(x_r)
&\!\!=\!\!&
\prod_{i<j}
\Exp{-\sum_{n>0}{1\/n}{1-t^n\/1-q^n}(1+p^n){x_j^n\/x_i^n}}
\prod_{i=1}^r x_a^{2\beta(r-i)}
:\prod_{i=1}^r S^a_+(x_i):\cr
&\!\!=\!\!&
\Delta(x) C(x)\prod_{i=1}^r x_i^{(r-1)\beta}:\prod_{i=1}^r S^a_+(x_i):,
\eea
and
\be
:\prod_{a=1}^{N-1}\prod_{i=1}^{r_a} S^a_+(x_i):|\widetilde\a_{r,s}\>
=\prod_{a=1}^{N-1}\prod_{i=1}^{r_a}
(x^a_i)^{(1-r_a+r_{a+1})\beta-(1+s_a)}(S^a_+(x_i))_-   \cdot|\a_{r,s}\>,
\ee
we obtain the proposition. \qed

Note that $C(x)$ is a pseudo-constant under the $q$-shift, {\it i.e.,}
$q^{D_{x_i}}C(x)=C(x)$. 
The expression in \eq{e:weightalpha}
is the same as that of $q=1$ case \cite{rAMOS}.


Remark that the singular vectors are also realized
by using the other screening currents $S_-^i(x)$
by the replacing $t$ with $q^{-1}$ and $\rb$ with $-1/\rb$ in
\eq{e:singular}, that is to say:
\bea
\&\hspace{-5mm}
|\chi_{r,s}^-\> =
\oint\prod_{a=1}^{N-1}\prod_{j=1}^{r_a} {dx^a_j}\cdot
S^1_-(x^1_1)\cdots S^1_-(x^1_{r_1}) \cdots
S^{N-1}_-(x^{N-1}_1) \cdots S^{N-1}_-(x^{N-1}_{r_{N-1}})
|\widetilde\a_{r,s}^-\>\cr
&\!\!=\!\!&
\oint\prod_{a=1}^{N-1} \prod_{j=1}^{r_a} {dx^a_j\/x^a_j}\cdot
\prod_{a=1}^{N-1} \Pi_-\(\overline{x^a},x^{a+1}\) \Delta_-(x^a) C_-(x^a)
\prod_{j=1}^{r_a} (x^a_j)^{-s_a} (S^a_-(x^a_j))_-\cdot|\a_{r,s}^-\>,\cr
&&\label{e:singularMinus}
\eea
where $\widetilde\a_{r,s}^-$, $\a_{r,s}^-$, $\Pi_-$, $\Delta_-$ and $C_-$
are obtained from those without $-$ suffix
by the replacing $t$ with $q$ and $\rb$ with $-1/\rb$.
And $(S^a_-(z))_-$ is the negative mode part of $S^a_-(z)$.


\section{Vertex operator of fundamental representation}


Now we introduce a vertex operator.
Let $V(z)$ be the vertex operator defined as
\be
V(z) \equiv
\,:\!\Exp{ -\sum_{n\neq0}{{\h 1n}\/1-q^n} p^{-{n\/2}}z^{-n} }\!:\,
	e^{-\rb{\Qh 1}} z^{-\rb{\h 10}}.
\label{e:vertex}
\ee
When $q=1$,
this $V(z)$ coincides with the vertex operator of fundamental representation.
Note that the fundamental vertex $\La_1(z)$ can be realized
by $V(z)$ as
\be
\La_1(zp^{1\/2}) = \,:\!V(zq^{-1}) V^{-1}(z)\!:\, p^{N-1\/2}.
\ee
Hence, this vertex operator $V(z)$ can be considered as
one of a building block of the {\qW} generators.
%
%
We have

\proclaim Proposition.
The vertex operator $V(w)$ enjoys the following Miura-like relation:
\bea
\&\hspace{-5mm}
:\!\(p^{D_z} - g^L\({w\/z       }\) \La_1(z       )\)\cdots
   \(p^{D_z} - g^L\({w\/zp^{1-N}}\) \La_N(zp^{1-N})\)\!:
V(w)\cr
\&
- V(w)
:\!\(p^{D_z} - \La_1(z       ) g^R\({z       \/w}\)\)\cdots
   \(p^{D_z} - \La_N(zp^{1-N}) g^R\({zp^{1-N}\/w}\)\)\!:
\cr
\& \hspace{5mm}
= p^{N-1\/2}(1-t^{-1})\ddelta\({w\/z}p^{1\/2}\)
:V(wq^{-1}) \(p^{D_z} - \La_2(zp^{-1})\)\cdots\(p^{D_z} - \La_N(zp^{1-N})\):,
\nonumber
\eea
 and
\bea
g^L(x) &\!\!=\!\!&
\Exp{\sum_{n>0} {1\/n}(1-t^{ n}){1-p^n   \/1-p^{ Nn}}p^{ {n\/2}}x^n }
t^{-{1\/N}},\cr
g^R(x) &\!\!=\!\!&
\Exp{\sum_{n>0} {1\/n}(1-t^{-n}){1-p^{-n}\/1-p^{-Nn}}p^{-{n\/2}}x^n }.
\eea

\proof
The fundamental relation is
\be
g^L\({w\/z}\) \La_i(z) V(w) - V(w) \La_i(z) g^R\({z\/w}\) =
p^{{N-1\/2}} (t^{-1}-1) \ddelta_{i,1} \ddelta\({w\/z}p^{1\/2}\) V(wq^{-1}),
\ee
{\it i.e.,}
\bea
\(p^{D_z} - g^L\({w\/z}\) \La_i(z)\)V(w)
&\!\!=\!\!&
V(w)\(p^{D_z} - \La_i(z) g^R\({z\/w}\)\)\cr
&\!\!+\!\!&
p^{N-1\/2}(1-t^{-1})\ddelta_{i,1}\ddelta\({w\/z}p^{1\/2}\) V(wq^{-1}),
\label{e:LaV}
\eea
here we use
$:\!\La_1(wp^{1\/2}) V(w)\!:\, = V(wq^{-1}) p^{{N-1\/2}}$.
By using this relation \eq{e:LaV} and
$V(w) \La_i(z) g^R\({z/w}\) = :\!V(w)\La_i(z)\!:$,
we obtain the proposition.
\qed


For example, when $N=3$,
the relation between the vertex operator $V(w)$ and the {\qW} generators is
\bea
\&
g^L\({w\/z}\) W^1(z) V(w) - V(w) W^1(z) g^R\({z\/w}\) 
=p (t^{-1}-1)\ddelta\({w\/z}p^{1\/2}\) V(wq^{-1}),\cr
\&
g^L\({w\/z}\) g^L\({w\/z}p\) W^2(zp^{-{1\/2}}) V(w) -
V(w) W^2(zp^{-{1\/2}}) g^R\({z\/w}\) g^R\({z\/w}p^{-1}\) \cr
\&\hspace{5mm}
=p (t^{-1}-1)\ddelta\({w\/z}p^{1\/2}\)
\(\,:\!V(wq^{-1})\La_2(wp^{-{1\/2}})\!:\,+
  \,:\!V(wq^{-1})\La_3(wp^{-{1\/2}})\!:\,\).
\eea


\section{Macdonald polynomials}


Finally we present a relation with the Macdonald polynomials.
The excited states of trigonometric Ruijsenaars model are
called Macdonald symmetric functions $P_\la(z)$
and they are defined as follows:
\bea
&&\qquad
 H P_\la(z_1,\cdots,z_M) =\varepsilon_\la P_\la(z_1,\cdots,z_M),\cr
&&
 H = \sum_{i=1}^M \prod_{j\neq i}
{t z_i - z_j \/ z_i - z_j}
\cdot q^{D_{z_i}},\qquad
 \varepsilon_\la = \sum_{i=1}^M t^{M-i} q^{\la_i},
\label{e:macDef}
\eea
where 
the $\la = (\la_1\geq\la_2\geq\cdots\la_M\geq0)$ is a partition.

The Macdonald polynomials with general Young diagram $\la$
are realized as some kinds of correlation functions of
the screening currents and vertex operators of the {\qW} algebra
as follows:

\proclaim Theorem.
Macdonald polynomial $P_\la(z)$ with the Young diagram
$\la = \sum_{i=1}^{N-1} (s_i^{r_i})$, $r_i\geq r_{i+1}$
is written as
\be
P_\la\(z_1,\cdots,z_M\)\propto
\<\a_{r,s}|\Exp{-\sum_{n>0}{{\h 1n}\/1-q^n}\sum_{i=1}^Mz_i^n}|\chi_{r,s}\>.
\ee
Here $|\chi_{r,s}\>$ is a singular vector in \eq{e:singular}.

Note that the operator part of the above equation
is the positive mode part of the product of the vertex operators \eq{e:vertex}.
The Young diagram is as follows:

\generalYoung

\proof
First we have
\be
\Exp{ -\sum_{n>0}{{\h 1n}\/1-q^n}\sum_{i=1}^M z_i^n} S^a_+(w) =
\Pi\(z,px^1\)^{\ddelta_{a,1}}
S^a_+(w) \Exp{ -\sum_{n>0}{{\h 1n}\/1-q^n}\sum_{i=1}^M z_i^n}.
\ee
By \eq{e:singular},
the right hand side of the equation of this theorem is
\be
\oint\prod_{a=1}^{N-1} \prod_{j=1}^{r_a} {dx^a_j\/x^a_j}\cdot
\Pi\(z,px^1\)
\prod_{a=1}^{N-1} \Pi\(\overline{x^a},px^{a+1}\) \Delta(x^a) C(x^a)
\prod_{j=1}^{r_a} (x^a_j)^{-s_a},
\label{e:macOPE}
\ee
If we replace $x^a$ with $(p^a x^a)^{-1}$ in \eq{e:macOPE},
then the integrand coincides with that of the integral formula for
Macdonald  polynomials in Ref. \cite{rAOS}
except for the $C(x)$ parts.
For the integral representation of the Macdonald polynomial,
we need only the property with respect to a $q$-shift.
Since this $C(x)$ is a pseudo-constant under it, {\it i.e.,}
$q^{D_{x_i}}C(x)=C(x)$, 
they are integral representations of the Macdonald polynomial
(see appendix B).
\qed


Remark that the Macdonald polynomials with the dual Young diagram
$\la'= \(r_1^{s_1},r_2^{s_2},\cdots,r_{N-1}^{s_{N-1}}\)$
are realized by using the other screening currents $S_-^i(x)$
with $|\chi_{r,s}^-\>$ in \eq{e:singularMinus} as
\be
P_{\la'}\(-z\)\propto
\<\a_{r,s}^-|\Exp{-\sum_{n>0}{{\h 1n}\/1-q^n}\sum_{i=1}^Mz_i^n}|\chi_{r,s}^-\>.
\ee


\section{Conclusion and discussion}

\def\FeiginFrenkel{
After finishing of this work,
we received the preprint
{\it ``Quantum $\cW$-algebras and elliptic algebras''}
by B.~Feigin and E.~Frenkel (q-alg/9508009).
They discuss similar things with ours.
Although the algebra of screening currents is considered there,
the normal ordering of $q$-$\cW$ generators and
the relation with the Macdonald polynomial are not given.
}

We have derived a quantum $\cW_N$ algebra
whose some kinds of correlation functions are the Macdonald polynomials.
\footnote\FeiginFrenkel

Jack polynomials are realized in the following two ways
(see also \cite{rLV}):
one is some kinds of correlation function of $\cW_N$ algebra
\cite{rMY,rAMOS},
the other is suitable combinations of
correlation functions of $\widehat{sl_N}$ algebra \cite{rMC}.
The relations between Macdonald polynomials,
the {\qW} algebra and the $U_q(\widehat{sl_N})$ algebra
are interesting.

In the classical limit $\hbar\rightarrow 0$ with $q\equiv e^\hbar$,
$q$-Miura transformation \eq{e:qMiura} reduces to
the classical one.
Since the right hand side of it is order $\hbar^N$,
the left hand side must be the same order.
To do so, $\hbar$ expansion of the {\qW} generators must be nontrivial.
Moreover, the classical generators are obtained as a linear
combination of the {\qW} generators.

\vskip5mm
\noindent{\bf Acknowledgments:}

\noindent
We would like to thank
B.~Feigin, E.~Frenkel and Y.~Matsuo 
for valuable discussions.
S.O. would like to thank members of YITP for their hospitality.
This work is supported in part by Grant-in-Aid for Scientific
Research from Ministry of Science and Culture. 


\section*{Appendix A: Quantum Virasoro algebra}

\def\Lukyanov{
The same operator with $S^1_+(z)$ was considered in \cite{rPL}.
}

In this appendix, we give an example when $N=2$,
{\it i.e.,} {\qVir} in \cite{rSKAO}.
The fundamental bosons ${\h 1n}$ and ${\Qh 1}$ satisfy
\be
[{\h 1n},{\h 1m}] = -{1\/n}{(1-q^n)(1-t^{-n})\/1+p^n}\ddelta_{n+m,0},\qquad
[{\h 10},{\Qh 1}] = {1\/2}.
\ee
The root bosons are
${\al 1n} = (1+p^{-n}) {\h 1n}$ and ${\Qal 1} = 2 {\Qh 1}$.

The $q$-Virasoro generator $W^1(z)$,
the screening currents $S^1_\pm(z)$ and
the vertex operator $V(z)$ are now\footnote\Lukyanov
\bea
W^1(z) &\!\!=\!\!&
:\Exp{ \sum_{n\neq 0}{\h 1n}      z^{-n}}:q^{ \rb{\h 10}} p^{  1\/2 }+
:\Exp{-\sum_{n\neq 0}{\h 1n}p^{-n}z^{-n}}:q^{-\rb{\h 10}} p^{-{1\/2}},\cr
S^1_\pm(z) &\!\!=\!\!&
:\Exp{\pm\sum_{n\neq0}{1+p^{-n}\/1-r_\pm^n}{\h 1n} z^{-n}}:
e^{\pm2\rb^{\pm1}{\Qh 1}} z^{\pm2\rb^{\pm1}{\h 10}},\quad
r_+ = q,\quad r_- = t,\cr
%
V(z) &\!\!=\!\!&
\,:\Exp{ -\sum_{n\neq0}{{\h 1n}\/1-q^n} p^{-{n\/2}}z^{-n} }:\,
	e^{-\rb{\Qh 1}} z^{-\rb{\h 10}}.
\eea

The relations of them are
\bea
f^{11}\({w\/z}\) W^1(z) W^1(w) &\!\!-\!\!& W^1(w) W^1(z) f^{11}\({z\/w}\) \cr
&\!\!=\!\!& -{(1-q)(1-t^{-1})\/1-p}
\left\{\ddelta\({w\/z}p\)-\ddelta\({w\/z}p^{-1}\)\right\},\cr
f^{11}(x) &\!\!=\!\!& \Exp{\sum_{n>0}{1\/n}{(1-q^n)(1-t^{-n})\/1+p^n}x^n},
\eea
\bea
\[\,W^1(z), S^1_\pm(w)\]
&\!\!=\!\!& -(1-q^{\pm1})(1-t^{\mp1})
{d\/d_{r_\pm}w} 
\left\{w\ddelta\({w\/z}\)A^1_\pm(w)\right\},\cr
A^1_\pm(w) &\!\!=\!\!&
:\Exp{\sum_{n\neq0}{1+r_\mp^{\pm n}\/1-r_\pm^{\pm n}}{\h 1n}w^{-n} }:\,
e^{\pm2\rb^{\pm1}{\Qh 1}}w^{\pm2\rb^{\pm1}{\h 10}}
q^{\mp\rb{\h 10}}p^{\mp{1\/2}},\nonumber
\eea
\bea
g^L\({w\/z}\) W^1(z) V(w) &\!\!-\!\!& V(w) W^1(z) g^R\({z\/w}\)
= p^{1\/2}(t^{-1}-1) \ddelta\({w\/z}p^{1\/2}\) V(wq^{-1}),\cr
g^{L\atop R}(x) &\!\!=\!\!&
\Exp{\sum_{n>0} {1\/n}{1-t^{\pm n}\/1+p^{\pm n}}p^{\pm{n\/2}}x^n }
t^{-{1\pm1\/4}}.
\eea

For non-negative integers $s$ and $r\geq0$,
the singular vectors $|\chi_{rs}\>\in\cF_{\a_{rs}}$ are
\bea
|\chi_{r,s}\> &\!\!=\!\!&
\oint\prod_{j=1}^{r} {dx_j}\cdot
S^1_+(x_1)\cdots S^1_+(x_{r}) |\a_{-r,s}\>\cr
&\!\!=\!\!&
\oint\prod_{j=1}^{r} {dx_j\/x_j}\cdot
\Delta(x) C(x)\prod_{j=1}^{r} (x_j)^{-s} (S_+(x_j))_-\cdot|\a_{r,s}\>,
\eea
with
$\a_{r,s}^1 =\rb(1+r)\rbi(1+s)$.
$\Delta(x)$ and $C(x)$ are the same as \eq{e:Delta}.


\section*{Appendix B: Integral formula for the Macdonald polynomials}


Finally, we recapitulate the integral representation
of the Macdonald polynomials \cite{rAOS}
(\cite{rMY2,rAMOS} in the $q=1$ case).
Let us denote the Macdonald polynomial defined by \eq{e:macDef} as
$P_\la(z;q,t)$ or $P_\la(z_1,\cdots,z_M;q,t)$.

\proclaim Proposition. 
The Macdonald polynomials with the Young diagram
$\la = \sum_{i=1}^{N-1} \(s_i^{r_i}\)$ or with its dual
$\la'= \(r_1^{s_1},r_2^{s_2},\cdots,r_{N-1}^{s_{N-1}}\)$ are
realized as follows:
\bea
P_\la(z;q,t)&\!\!\propto\!\!&
\oint\prod_{a=1}^{N-1} \prod_{j=1}^{r_a} {dx^a_j\/x^a_j}\cdot
\Pi\(z,\overline{x^1}\)
\prod_{a=1}^{N-1} \Pi\(x^a,\overline{x^{a+1}}\) \Delta(x^a) C(x^a)
\prod_{j=1}^{r_a} (x^a_j)^{s_a},\cr
P_{\la'}(z;t,q)&\!\!\propto\!\!&
\oint\prod_{a=1}^{N-1} \prod_{j=1}^{r_a} {dx^a_j\/x^a_j}\cdot
\widetilde\Pi\(z,\overline{x^1}\)
\prod_{a=1}^{N-1} \Pi\(x^a,\overline{x^{a+1}}\) \Delta(x^a) C(x^a)
\prod_{j=1}^{r_a} (x^a_j)^{s_a},\nonumber
\eea
with an arbitrary pseudo-constant $C(x)$ such that $q^{D_{x_i}}C(x)=C(x)$.
Here $\widetilde\Pi(x,y)\equiv \prod_{ij}(1+x_i y_j)$.
$\Pi$ and $\Delta$ are in \eq{e:Delta}.

\proof
This proposition is proved by using
two transformations in the following lemmas iteratively.
The first transformation adds a rectangle to the Young diagram
and the second one increases the number of variables.
\qed


\proclaim Lemma 1. Galilean transformation.
$($eq.\ $(VI.4.17)$ in $\cite{rM})$
\be
P_{\la+(s^r)}(x_1,\cdots,x_r) = P_\la(x_1,\cdots,x_r)\prod_{i=1}^r x_i^s.
\ee

This transformation adds a rectangle Young diagram to the original one:
$$
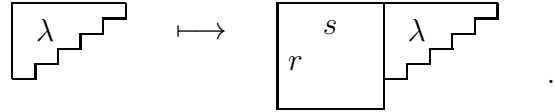
\Galilei
$$


\proclaim Lemma 2. Particle number changing transformation.
\bea
P_\la(x_1,\cdots,x_N;q,t) &\!\!\propto\!\!&
\oint \prod_{j=1}^M {dy_j\/y_j}
\Pi(x,\overline y)\Delta(y) C(y) P_\la(y_1,\cdots,y_M;q,t),\cr
P_{\la'}(x_1,\cdots,x_N;t,q) &\!\!\propto\!\!&
\oint \prod_{j=1}^M {dy_j\/y_j}
\widetilde\Pi(x,\overline y)\Delta(y) C(y) P_\la(y_1,\cdots,y_M;q,t),\nonumber
\eea
here $C(y)$ is an arbitrary pseudo-constant $q^{D_{y_i}}C(y)=C(y)$
and $\la'$ is a dual Young diagram of $\la$.

\proof
Let us define scalar products $\<*,*\>$ and the another one $\<*,*\>'_N$
as follows:
\bea
\<f,g\>    &\!\!\equiv\!\!&
\oint\prod_{n>0}   {dp_n\/2\pi i p_n}\,{f(\overline p)}\,g(p),\cr
\<f,g\>'_N &\!\!\equiv\!\!& {1\/N!}
\oint\prod_{j=1}^N {dx_j\/2\pi i x_j}\Delta(x)\,{f(\overline x)}\,g(x),
\eea
for the symmetric functions $f$ and $g$
with $p_n\equiv\sum_{i=1}^N x_i^N$,
$\overline{p_n}\equiv n{1-q^n\/1-t^n}{\partial\/\partial p_n}$ and
$\overline{x_j}\equiv{1/x_j}$ .
Here we must treat the power-sums $p_n$ as formally independent variables,
{\it i.e.},
${\partial\/\partial p_n}\, p_m = \delta_{n,m}$ for all $n,m>0$.
Then  (eq.\ (VI.4.13) and (VI.5.4) in \cite{rM})
\bea
\Pi(x,y) &\!\!=\!\!&
\sum_\la P_\la(x;q,t) P_\la(y;q,t) \<P_\la,P_\la\>^{-1},\cr
\widetilde\Pi(x,y) &\!\!=\!\!&
\sum_\la P_\la(x;q,t) P_{\la'}(y;t,q).
\label{e:completeness}\eea
Since the Macdonald operator is self-adjoint
for the another scalar product $\<*,*\>'_N$, that is to say
$\<H\,f,g\>'_N = \<f,H\,g\>'_N$ (eq.\ (VI.9.4) in \cite{rM}),
the Macdonald polynomials are orthogonal for this product
$\<P_\la,C\,P_\mu\>'_N \propto \delta_{\la,\mu}$ 
with an arbitrary pseudo-constant $C$.
The proposition follows from the completeness \eq{e:completeness} and
the orthogonality of $P_\la$'s.
\qed


Remark that the above lemma 2 is also proved directly by using
the power-sum representation of the Macdonald operator \cite{rAMOS}.
Since that is also important to analyze the algebraic properties of
the Macdonald polynomials, we review it here.

\proclaim Proposition.
Macdonald operator $H(x_1,\cdots,x_N)$ are written
by the power sums $p_n \equiv \sum_{i=1}^N x_i^n$ as follows:
\be
H = {t^N\/t-1}\oint{d\xi\/2\pi i \xi}
\Exp{\sum_{n>0}{1-t^{-n}\/n}p_n \xi^n}
\Exp{\sum_{n>0}(q^n-1){\partial\/\partial p_n} \xi^{-n}}
-{1\/t-1}.
\ee

\proof
Since $q^{D_{x_i}} p_n = \((q^n-1)x_i^n+p_n\) q^{D_{x_i}}$,
we have 
\be
q^{D_{x_i}} =\,:\!\Exp{\sum_{n>0}(q^n-1)x_i^n{\partial\/\partial p_n}}\!:\,
= \oint{d\xi\/2\pi i \xi} \sum_{n\geq0} x_i^n \xi^n \cdot
\Exp{\sum_{n>0}(q^n-1){\partial\/\partial p_n} \xi^{-n}},
\ee
here $:*:$ stands for the normal ordering such that
the differential operators ${\partial\/\partial p_n}$ are in the right.
It follows from eq.\ (III.2.9) and (III.2.10) in \cite{rM} that
\be
\sum_i\prod_{j\neq i} {tx_i-x_j\/x_i-x_j}
\sum_{n\geq 0}x_i^n \xi^n
=
{t^N\/t-1}\Exp{\sum_{n>0}{1-t^{-n}\/n}p_n \xi^n}
-{1\/t-1}.
\ee
This gives us the proposition. \qed

Let
$\widetilde H_N(x_1,\cdots,x_N) \equiv t^{-N}\((t-1)H(x_1,\cdots,x_N)+1\)$,
then
\be
\widetilde H_N(x_1,\cdots,x_N) \Pi(x,y) =
\widetilde H_M(y_1,\cdots,y_M) \Pi(x,y).
\ee
With the self-adjointness of $H$ for the another scalar product,
we obtain the lemma 2 again.





\begin{thebibliography}{99}


\bibitem{rAMOS}
{H.~Awata, Y.~Matsuo, S.~Odake and J.~Shiraishi,
{\sl Phys. Lett.} {\bf B347} (1995) 49-55; 
{\it Excited States of Calogero-Sutherland Model and
	Singular Vectors of the $W_N$ Algebra}, hep-th/9503043,
to appear in {\sl Nucl. Phys.} {\bf B};
{\it A Note on Calogero-Sutherland Model, $W_n$ Singular Vectors
   and Generalized Matrix Models}, hep-th/9503028,
  to appear in  Soryushiron kenkyu (Kyoto)}

\bibitem{rAOS}
{H.~Awata, S.~Odake and J.~Shiraishi,
{\it Integral Representations of the Macdonald Symmetric Functions},
q-alg/9506006}

\bibitem{rFL}
{V. Fateev and S. Lukyanov, Int. J. Mod. Phys. {\bf A3} (1988) 507--520}

\bibitem{rFR}
{E.~Frenkel and N.~Reshetikhin,
{\em Quantum Affine Algebras and Deformations
of The Virasoro and $\cal W$-Algebra}, q-alg/9505025}

\bibitem{rLV}{L. Lapointe and L. Vinet
{\it Exact Operator Solution of the Calogero-Sutherland model},
hep-th/9507073}

\bibitem{rM}{I.G.~Macdonald
{\sl ``Symmetric Functions and Hall Functions''} (2nd ed.),
Oxford University Press 1995}

\bibitem{rMC}{A.~Matsuo, {\sl Invent. Math.} {\bf 110} (1992) 95-121;\hb
  I.~Cherednik,
  {\it Integration of quantum many-body problems
   by affine Knizhnik-Zamolodchikov equations},
 preprint, RISM-776 (1991)}

\bibitem{rMY}
{K.~Mimachi and Y.~Yamada,
  {\it Singlar vectors of the Virasoro algebra in terms of
  Jack symmetric polynomials},
  Kyushu univ. preprint (November 1994), 
 to appear in {\sl Commun. Math. Phys.}}

\bibitem{rMY2}
{K.~Mimachi and Y.~Yamada,
 Talk at the workshop ``Hypergeometric Functions''
at RIMS, 1994 Dec. 12--15}

\bibitem{rPL}{Y. Pugai and S. Lukyanov,
{\it Bosonization of ZF Algebras: Direction Toward Deformed Virasoro Algebra},
hepth/9412128}

\bibitem{rR}{S.N.M.~Ruijsenaars,
  {\sl Comm. Math. Phys.} {\bf 110} (1987) 191-213}

\bibitem{rSKAO}
{J. Shiraishi, H. Kubo, H. Awata and S. Odake,
{\em A Quantum Deformation of the Virasoro Algebra and
the Macdonald Symmetric Functions}, q-alg/9507034,
to appear in {\sl Lett. Math. Phys.}}

\bibitem{rSt}{R. Stanley, {\sl Adv. Math.} {\bf 77} (1989) 76--115}

\bibitem{rSu}
{B.~Sutherland,
{\sl Phys. Rev.} {\bf A4} (1971) 2019--2021; {\bf A5} (1972) 1372--1376}


\end{thebibliography}
\end{document}